\documentclass[a4paper,11pt]{amsart}
\usepackage{amscd,amssymb}
\usepackage[all]{xy}

\newcommand{\hhdim}{\operatorname{hhdim}}
\newcommand{\chdim}{\operatorname{chdim}}

\newcommand{\HC}{\operatorname{HC}}
\newcommand{\HH}{\operatorname{HH}}
\newcommand{\Tor}{\operatorname{Tor}}
\newcommand{\Ext}{\operatorname{Ext}}

\newcommand{\Z}{\mathbb Z}
\newcommand{\e}{\operatorname{e}}
\newcommand{\op}{\operatorname{op}}
\newcommand{\gr}{\operatorname{gr}}
\newcommand{\az}{\operatorname{\mathfrak{a}}\nolimits}

\newtheorem{lem}{Lemma}[section]
\newtheorem{prop}[lem]{Proposition}
\newtheorem{cor}[lem]{Corollary}
\newtheorem{thm}[lem]{Theorem}

\theoremstyle{definition}

\newtheorem{example}[lem]{Example}

\newtheorem*{conj}{Conjecture}

\begin{document}

\title{Hochschild homology and global dimension}

\author{Petter Andreas Bergh \& Dag Madsen}
\address{Institutt for matematiske fag\\ NTNU\\ 7491 Trondheim\\ Norway}
\email{bergh@math.ntnu.no} \email{dagma@math.ntnu.no}

\thanks{The authors were supported by NFR Storforsk grant no.\
167130}



\maketitle

\begin{abstract}
We prove that for certain classes of graded algebras (Koszul, local,
cellular), infinite global dimension implies that Hochschild
homology does not vanish in high degrees, provided the
characteristic of the ground field is zero. Our proof uses Igusa's
formula relating the Euler characteristic of relative cyclic
homology to the graded Cartan determinant.
\end{abstract}

\section{Introduction}

The homological properties of a finite dimensional algebra are
closely related to the behavior of the algebra as a bimodule over
itself. For example, if the algebra has finite projective dimension
as a bimodule, then its global dimension is also finite. The
converse holds if the algebra modulo its Jacobson radical is
separable over the ground field, something which automatically
happens when the field is algebraically closed.

In particular, if a finite dimensional algebra over an algebraically
closed field has finite global dimension, then all its higher
Hochschild cohomolgy groups vanish. In \cite{Happel}, following this
easy observation, Happel remarked that ``the converse seems to be
not known'', thus giving birth to what subsequently became known as
``Happel's question": if all the higher Hochschild cohomolgy groups
of a finite dimensional algebra vanish, then is the algebra of
finite global dimension? As shown in \cite{AvramovIyengar}, the
answer is yes when the algebra is commutative. However, it was shown
in \cite{Buchweitz} that the answer in general is no. Namely, given
a field $k$ and a nonzero element $q \in k$ which is not a root of
unity, then the total Hochschild cohomology of the four-dimensional
algebra
$$k \langle X,Y \rangle / (X^2, XY-qYX, Y^2)$$
is five. In particular, all the higher Hochschild cohomology groups
of this algebra vanish, whereas the algebra, being selfinjective,
clearly does not have finite global dimension.

As shown by Han in \cite{Han}, the total Hochschild homology of the
above algebra is infinite dimensional. Han then conjectured that the
homology version of Happel's question would always hold, namely that
a finite dimensional algebra whose higher Hochschild homology groups
vanish must be of finite global dimension. In the same paper, he
showed that the conjecture holds for monomial algebras. Moreover, as
in the cohomology case, the conjecture holds if the algebra is
commutative, by \cite{AvramovPoirrier} (for finitely generated but
not necessarily finite dimensional algebras, see also
\cite{Poirrier}).

In this paper, we show that Han's conjecture holds for graded local
algebras, Koszul algebras and graded cellular algebras, provided the
characteristic of the ground field is zero. We do this by exploiting
some particular properties of the graded Cartan matrix and the
logarithm of its determinant, concepts extensively studied in
\cite{Igu}.

\section{Hochschild homology and cyclic homology}

Throughout this paper, we fix a field $k$, not necessarily
algebraically closed. Let $A$ be a finite dimensional $k$-algebra,
and denote by $A^{\e}$ the enveloping algebra $A \otimes_k A^{\op}$
of $A$. The \emph{Hochschild homology} of $A$, denoted $\HH_* (A)$,
is defined by $\HH_* (A) = \Tor^{A^{\e}}_* (A,A)$. By definition, it
is obtained by taking any projective bimodule resolution of $A$,
applying $A \otimes_{A^{\e}} -$ and computing the homology of the
resulting complex. However, we shall explore one particular such
resolution, which eventually leads to the definition of cyclic
homology.

For each $n \ge 0$, denote by $A^{\otimes n}$ the $n$-fold tensor
product $A \otimes_k \cdots \otimes_k A$, in which there are $n$
copies of $A$. For $n \ge 2$, this is a projective bimodule. Define
a map
\begin{eqnarray*}
A^{\otimes (n+1)} & \xrightarrow{b'} & A^{\otimes n} \\
a_0 \otimes \cdots \otimes a_n & \mapsto & \sum_{i=0}^{n-1} (-1)^i
a_0 \otimes \cdots \otimes a_ia_{i+1} \otimes \cdots \otimes a_n,
\end{eqnarray*}
and consider the complex
$$\cdots \to A^{\otimes 4} \xrightarrow{b'} A^{\otimes 3}
\xrightarrow{b'} A^{\otimes 2}$$ in which $A^{\otimes 2}$ is in
degree zero. By \cite[\S IX.6]{Cartan} this complex is exact, and it
is therefore a projective bimodule resolution of $A$. This is the
\emph{standard resolution} (or \emph{Bar resolution}) of $A$.
Applying $A \otimes_{A^{\e}} -$ to this resolution, we obtain the
complex
$$\cdots \to A^{\otimes 3} \xrightarrow{b} A^{\otimes 2}
\xrightarrow{b} A,$$ in which the map $A^{\otimes (n+1)}
\xrightarrow{b} A^{\otimes n}$ is given by
\begin{eqnarray*}
a_0 \otimes \cdots \otimes a_n & \mapsto & \sum_{i=0}^{n-1} (-1)^i
a_0 \otimes \cdots \otimes a_ia_{i+1} \otimes \cdots \otimes a_n \\
& & + (-1)^n a_na_0 \otimes a_1 \otimes \cdots \otimes a_{n-1}.
\end{eqnarray*}
By definition, the homology of this complex is the Hochschild
homology of our algebra $A$.

The standard resolution and the Hochschild complex are also the key
ingredients in the definition of cyclic homology. For each $n \ge
0$, define the map
\begin{eqnarray*}
A^{\otimes (n+1)} & \xrightarrow{t} & A^{\otimes (n+1)} \\
a_0 \otimes \cdots \otimes a_n & \mapsto & a_n \otimes a_0 \otimes
\cdots \otimes a_{n-1},
\end{eqnarray*}
and let $N = 1+t+ \cdots +t^n$ be the corresponding norm operator.
Then $(1-t)b' = b(1-t)$ and $b'N=Nb$ (cf.\ \cite[Lemma
2.1.1]{Loday}), and so
$$\xymatrix{
\vdots \ar[d] & \vdots \ar[d] & \vdots \ar[d] \\
A^{\otimes 3} \ar[d]^b & A^{\otimes 3} \ar[d]^{-b'} \ar[l]_{1-t} &
A^{\otimes 3} \ar[d]^{b} \ar[l]_{N} & \cdots \ar[l]_{1-t} \\
A^{\otimes 2} \ar[d]^b & A^{\otimes 2} \ar[d]^{-b'} \ar[l]_{1-t} &
A^{\otimes 2} \ar[d]^{b} \ar[l]_{N} & \cdots \ar[l]_{1-t} \\
A & A \ar[l]_{1-t} & A \ar[l]_{N} & \cdots \ar[l]_{1-t} }$$ is a
first quadrant double complex (in which the lower left $A$ has
degree $(0,0)$). The \emph{cyclic homology} of $A$, denoted $\HC_*
(A)$, is the homology of the resulting total complex. It is closely
linked to the Hochschild cohomology of $A$ via the well known long
exact sequence
$$\cdots \to \HH_n (A) \xrightarrow{I} \HC_n (A) \xrightarrow{S}
\HC_{n-2} (A) \xrightarrow{B} \HH_{n-1} (A) \xrightarrow{I} \cdots$$
due to Connes, the \emph{SBI sequence} or \emph{Connes' exact
sequence}.

Both Hochschild and cyclic homology are functorial, since an algebra
homomorphism induces a map between the corresponding Hochschild
complexes and cyclic double complexes. In particular, if $\az$ is a
twosided ideal of $A$, then the surjection $A \to A / \az$ induces a
surjective map of cyclic double complexes. The kernel of this map is
a double complex, and the homology of its total complex, denoted
$\HC_* (A, \az )$, is the \emph{relative cyclic homology} of $A$
with respect to $\az$.

Suppose the algebra $A$ is graded, say $A = A_0 \oplus \cdots \oplus
A_s$. Then its grading induces an internal grading
\begin{eqnarray*}
\HH_* (A) & = & \oplus_i \HH_*^i (A) \\
\HC_* (A) & = & \oplus_i \HC_*^i (A) \\
\HC_* (A, \az ) & = & \oplus_i \HC_*^i (A, \az)
\end{eqnarray*}
on both Hochschild and (relative) cyclic homology. Consider the SBI
sequence relating the Hochschild and cyclic homology of $A$. By a
theorem of Goodwillie (cf.\ \cite[Corollary II.4.6]{Good} or
\cite[Theorem 9.9.1]{Weibel}), the image of the map $\HC_n^i (A)
\xrightarrow{S} \HC_{n-2}^i (A)$ is annihilated by $i$ (where $i$ is
viewed as an element in $k$) for every $i >0$. In particular, if the
characteristic of $k$ is zero, then we obtain a short exact sequence
\begin{equation*}\label{ES}
0 \to \HC_{n-1}^i(A) \to \HH_n^i(A) \to \HC_n^i(A) \to 0
\tag{$\dagger$}
\end{equation*}
for each $i>0$. It follows from this sequence that if $\HC_n^i (A)$
is nonzero for some $i>0$, then so is $\HH_{n+1}^i (A)$. Next,
suppose $A_0$ is a product of copies of the ground field $k$, say
$A_0 = k^{\times r}$, and let $J$ denote the radical $A_1 \oplus
\cdots \oplus A_s$ of $A$. If the characteristic of $k$ is zero,
then it folloes from \cite[Corollary 1.2]{Igu} that $\HC_n^m (A,J)
=0$ for any $m \ge 1$ and $n \ge m$. Therefore, in this case, the
\emph{Euler characteristic} of $\HC_*^m (A,J)$ is well defined and
given by
$$\chi \left ( \HC_*^m (A,J) \right ) =
\sum_{n=0}^{m-1} (-1)^i \dim_k \HC_n^m (A,J),$$ and we define the
\emph{graded Euler characteristic} by
$$\chi ( \HC_* (A,J))(x) \stackrel{\text{def}}{=}
\sum_{m=1}^{\infty} \chi ( \HC_*^m (A,J)) x^m.$$ The latter is a
power series with integer coefficients.

The aim of this paper is to establish results concerning the
non-vanishing of Hochschild homology for certain algebras. To
simplify the notation, we therefore define the following:
\begin{eqnarray*}
\hhdim A & \stackrel{\text{def}}{=} & \sup \{ n \in \mathbb{Z}
\mid \HH_n (A) \neq 0 \} \\
\chdim A & \stackrel{\text{def}}{=} & \sup \{ n \in \mathbb{Z}
\mid \HC_n (A) \neq 0 \} \\
\chdim (A, \az ) & \stackrel{\text{def}}{=} & \sup \{ n \in
\mathbb{Z} \mid \HC_n (A, \az ) \neq 0 \}.
\end{eqnarray*}
As mentioned in the introduction, we show in this paper that if the
characteristic of $k$ is zero and $A$ has infinite global dimension,
then $\hhdim A = \infty$ when $A$ is graded local, Koszul or graded
cellular. We end this section with the following result, which shows
that we only need to establish the non-vanishing of the relative
cyclic homology of $A$ with respect to the radical.

\begin{lem}\label{relvan}
Suppose $A$ is graded and that its degree zero part is a product of
copies of $k$. Furthermore, suppose the characteristic of $k$ is
zero. Then
$$\chdim (A,J) = \infty \hspace{.1cm} \Leftrightarrow \hspace{.1cm} \hhdim A =
\infty,$$ where $J$ is the radical of $A$.
\end{lem}

\begin{proof}
\sloppy By construction, there is a long exact sequence
$$\cdots \to \HC_{n+1} (A_0) \to \HC_n (A,J) \to \HC_n (A) \to \HC_n
(A_0) \to \cdots$$ relating relative and ordinary cyclic homology.
Since $A_0$ lives only in degree zero, its internal positive degree
cyclic homology vanishes, i.e.\ $\HC_*^m(A_0)=0$ for $m>0$.
Consequently, there is an isomorphism $\HC_n^m(A,J) \simeq
\HC_n^m(A)$ for every $n$ and any $m>0$. Suppose $\HC_n (A,J)$ is
nonzero for some $n$. Since $A_0$ is a product of copies of $k$, we
know from above that $\HC_n^i (A,J)=0$ for $i \le n$, and therefore
there must be an integer $m>n$ such that $\HC_n^m (A,J)$ is nonzero.
The isomorphism above then shows that $\HC_n^m(A)$ is nonzero, and
from the exact sequence (\ref{ES}) we see that the same holds for
$\HH_n^m(A)$. This shows the implication $\Rightarrow$. For the
reverse implication, note that for $n>0$, the group $\HH_n^0(A)$
vanishes. Hence if $\HH_n(A)$ is nonzero, then from the SBI sequence
we see that there is an $m>0$ such that either $\HC_n^m(A)$ or
$\HC_{n-1}^m(A)$ is nonzero. Then either $\HC_n^m(A,J)$ or
$\HC_{n-1}^m(A,J)$ must be nonzero.
\end{proof}

\section{The graded Cartan determinant}

In this section $A$ denotes a positively graded finite dimensional
$k$-algebra $A=A_0 \oplus A_1 \oplus \dots \oplus A_s$. We assume
$A_0 \simeq k \times \dots \times k=k^{\times r}$ as rings. Let
$1_A=e_1+ \ldots +e_r$ be the corresponding decomposition of the
identity.

In \cite{Igu}, the author presented a formula relating the graded
Euler characteristic $\chi (\HC_\ast(A,J))(x)$ of relative cyclic
homology to the so-called graded Cartan determinant of $A$. For $1
\leq l \leq s$, let $C^l$ be the $r \times r$ matrix with entries
$C^l_{i,j}=\dim_k e_j A_l e_i$. The \emph{graded Cartan matrix} of
$A$ is defined to be the $r \times r$ matrix
$$C_A(x)=C^0+C^1 x+C^2 x^2+ \ldots + C^s x^s$$ with entries in $\Z
[x]$. Its determinant $\det C_A(x)$ is the \emph{graded Cartan
determinant} of $A$. With our assumptions $C^0$ is the identity
matrix, and therefore $\det C_A(x)$ is a polynomial of degree $u
\leq s$ with integer coefficients and constant term $1$. The
logarithm of the determinant is then a power series that can be
defined according to the formula
$$\log \det C_A (x) = \sum_{m=1}^{\infty} (-1)^{m+1} \frac {(\det C_A (x)-1)^m} m.$$
Although this power series in general has rational coefficients, its
formal derivative $D_x (\log \det C_A (x))$ has integer
coefficients. We say that a power series is \emph{proper} if it has
infinitely many nonzero terms.

\begin{lem}\label{derlemma}
The power series $D_x (\log \det C_A (x))$ is proper and has integer
coefficients $\{b_i\}_{i \geq 0}$. The sequence $\{b_i\}_{i \geq 0}$
satisfies a linear recurrence relation of order $u$ with constant
integer coefficients.
\end{lem}

\begin{proof}
The chain rule gives $$D_x (\log \det C_A (x)) \cdot \det C_A
(x)=D_x(\det C_A (x)),$$ or alternatively $$D_x (\log \det C_A (x))
= (\det C_A (x))^{-1} \cdot D_x(\det C_A (x)).$$ Since the degree of
the polynomial $D_x(\det C_A (x))$ is strictly less than the degree
of $\det C_A (x)$, it follows from the first formula that $D_x (\log
\det C_A (x))$ must be a proper power series. Since $D_x(\det C_A
(x))$ has integer coefficients, it follows from the second formula
that $D_x (\log \det C_A (x))$ has integer coefficients.

If $\det C_A (x)=1+c_1 x+ \ldots + c_u x^u$, then for any $m\geq u$,
the first formula gives $b_m+c_1 b_{m-1}+ \ldots+c_u b_{m-u}=0$. The
last statement of the lemma follows.
\end{proof}

Next we state Igusa's formula as presented in \cite{Igu}. Recall
that the M\"{o}bius function $\mu$ is the multiplicative number
theoretic function defined by
$$\mu(n)=\left \{
\begin{array}{l l}
1 & \text{if $n=1$,}\\
(-1)^t & \text{if $n$ is a product of $t$ distinct primes,}\\
0 & \text{if $n$ has one or more repeated prime factors.}
\end{array} \right.$$

\begin{thm}\cite[Theorem 3.5]{Igu}\label{iguthm}
Let $A$ be a graded algebra over a field $k$ of characteristic zero,
and suppose $A_0$ is a product of copies of $k$. Then
\begin{itemize}
\item[(a)] $$\chi (\HC_\ast(A,J))(x)
=\sum_{m=1}^{\infty} \log \det C_A (x^m) \sum_{d \mid m} \frac {\mu
(d)} d,$$
\item[(b)] $$\log \det C_A (x)
=\sum_{m=1}^{\infty} \chi (\HC_\ast(A,J))(x^m) \sum_{d \mid m} \frac
{d \mu (d)} m,$$ where $\mu$ is the M\"{o}bius function.
\end{itemize}
\end{thm}

An immediate corollary of this theorem is that $\det C_A(x)=1$ if
and only if $\chi (\HC_\ast(A,J))(x)=0$. Our aim is to show that if
$\det C_A(x) \neq 1$, then $\chi (\HC_\ast(A,J))(x)$ is a proper
power series.

We would like to have a formula relating the coefficients of $$D_x
(\log \det C_A (x))=\sum_{i=1}^{\infty} b_i x^i$$ with the
coefficients of $$\chi (\HC_\ast(A,J))(x)=\sum_{i=1}^{\infty} a_i
x^i.$$ For this purpose we introduce the number theoretic function
$\theta$, defined as follows.
$$\theta (m) =\sum_{d \mid m} d \mu (d)=\prod_{\substack{p \mid m \\p \text{ prime}}} (1-p).$$
The function $\theta$ is multiplicative and has the property that if
$n \mid m$, then $\theta (n) \mid \theta (m)$. From Theorem
\ref{iguthm} (b) we get
$$D_x (\log \det C_A (x))=\sum_{m=1}^{\infty} [\sum_{d \mid m} a_d \cdot d \cdot \theta (\frac m d)] x^{m-1}.$$
For convenience let $$f(m)=\sum_{d \mid m} a_d \cdot d \cdot \theta
(\frac m d).$$ With this notation $f(m)=b_{m-1}$.

If $\chi (\HC_\ast(A,J))(x)$ has only finitely many non-zero
coefficients, in other words if $\chi (\HC_\ast(A,J))(x)$ is a
polynomial, then we get the following condition on the function $f$.

\begin{prop}\label{conprop}
If $\chi (\HC_\ast(A,J))(x)$ is a polynomial, then for every pair of
integers $s,t>0$ there are $s$ consecutive positive integers $N+1,
\ldots, N+s$ such that $2^t \mid f(N+i)$, $1 \leq i \leq s$.
\end{prop}

\begin{proof}
Let $v$ be the degree of $\chi (\HC_\ast(A,J))(x)$. Choose $st$
different odd prime numbers $p_{i,j}>v$, $1 \leq i \leq s, 1 \leq j
\leq t$. Let $n_i=\prod_{j=1}^t p_{i,j}$ for each $1 \leq i \leq s$.
The system of congruences
$$x \equiv -1 \mod n_1,$$
$$x \equiv -2 \mod n_2,$$
$$\vdots$$
$$x \equiv -s \mod n_s$$
has a unique solution modulo $n_1 n_2 \cdots n_s$. Let $N > 0$ be a
solution. Now for each $1 \leq i \leq s$, we have $n_i \mid N+i$.
Also $n_i \mid \frac {N+i} d$ whenever $d \leq v$. Since $a_d=0$ for
$d>v$, we have
$$f(N+i)=\sum_{\substack {d \mid (N+i)\\d \leq v}} a_d \cdot d \cdot \theta (\frac {N+i} d).$$
Since $\theta(n_i) \mid \theta (\frac {N+i} d)$ whenever $d \leq v$,
it follows that $\theta (n_i) \mid f(N+i)$. Since $\theta
(n_i)=\prod_{j=1}^t (1-p_{i,j})$, it follows that $2^t \mid \theta
(n_i) \mid f(N+i)$.
\end{proof}

If $\det C_A \neq 1$, then $f$ cannot satisfy this condition, and as
a consequence we get our theorem.

\begin{thm}\label{mainthm}
Let $A$ be a graded finite dimensional algebra over a field $k$ of
characteristic zero, and suppose $A_0$ is a product of copies of
$k$. If $\det C_A \neq 1$, then $\chi (\HC_\ast(A,J))(x)$ is a
proper power series.
\end{thm}

\begin{proof}
Suppose for contradiction that $\det C_A \neq 1$ and $\chi
(\HC_\ast(A,J))(x)$ is a polynomial. Since $\det C_A \neq 1$, it
follows that $\chi (\HC_\ast(A,J))(x) \neq 0$. Let $v$ be the degree
of $\chi (\HC_\ast(A,J))(x)$, and let $l$ be the degree of its
lowest degree nonzero term. Let $p>v$ be a prime number. Then
$f(lp^i)=a_l \cdot l \cdot \theta (p^i)=a_l \cdot l (1-p)=f(lp) \neq
0$ for any $i \geq 1$. Let $t \geq 1$ be the number such that
$2^{t-1} \mid f(lp)$ but $2^t \nmid f(lp)$.

Let $u>0$ be the degree of $\det C_A (x)$. By Proposition
\ref{conprop} there are $u$ consecutive positive integers $N+1,
\ldots, N+u$ such that $2^t \mid f(N+i)$, $1 \leq i \leq u$. So $D_x
(\log \det C_A (x))$ has $u$ consecutive coefficients
$b_N,\ldots,b_{N+u-1}$ which are divisible by $2^t$. But then by
Lemma \ref{derlemma} we have $2^t \mid b_m$ for any $m \geq N$, so
$2^t \mid f(m)$ for any $m \geq N+1$. This contradicts the fact that
$2^t \nmid f(lp^i)$ for any $i \geq 0$.
\end{proof}

For our investigations into the validity of Han's conjecture, the
following corollary is important.

\begin{cor}\label{determinant}
Let $A$ be a graded finite dimensional algebra over a field $k$ of
characteristic zero, and suppose $A_0$ is a product of copies of
$k$. Suppose $\det C_A \neq 1$. Then
\begin{itemize}
\item[(a)] $\chdim (A,J)=\infty$,
\item[(b)] $\hhdim A=\infty$.
\end{itemize}
\end{cor}

\begin{proof}
Parts (a) and (b) are equivalent by Lemma \ref{relvan}. Since $A$ is
finite dimensional, for each $n \geq 0$ the group $\HC_n (A,J)$ is
finite dimensional and therefore $\HC_n^m(A,J) \neq 0$ only for
finitely many internal degrees $m$. If $\chi (\HC_\ast(A,J))(x)$ is
a proper power series, then $\HC_\ast^m (A,J) \neq 0$ for infinitely
many internal degrees $m$, and this is only possible if $\chdim
(A,J)=\infty$. So it follows from Theorem \ref{mainthm} that if
$\det C_A \neq 1$, then $\chdim (A,J)=\infty$.
\end{proof}

\section{Non-vanishing of Hochschild homology in high degrees}

In this section we apply Corollary \ref{determinant} to various
classes of graded algebras (Koszul, local, cellular), and prove that
Han's conjecture \cite{Han} holds for these classes. Han's
conjecture can be stated as follows.

\begin{conj}[Han] Let $A$ be a finite dimensional algebra over a field.
If $A$ has infinite global dimension, then $\hhdim A = \infty$.
\end{conj}

Let $A$ be a graded finite dimensional $k$-algebra, and suppose $A_0
\simeq k^{\times r}$. The matrix $C_A(x)$ is invertible as a matrix
over $\Z[x]$ if and only if $\det C_A(x)=1$. On the other hand,
since $\det C_A(x)$ has constant term $1$, the matrix $C_A(x)$ is
always invertible when considered as a matrix over $\Z [[x]]$. In
\cite{Wil} we find the following formula for the entries in the
inverse matrix $C_A^{-1}(x)$. Here $\Ext_{\gr A}$ denotes graded
extensions and $X[u]$ denotes the $u$th shift of the graded module
$X$. For each $1 \leq i \leq r$, we let $S_i$ denote the degree zero
simple module $S_i=A e_i/J A e_i$.

\begin{thm}\cite[Theorem 1.5]{Wil}
The entry $c_{ij}$ in $C_A^{-1}(x)$ is given by
$$c_{ij} = \sum_{u \geq 0} \sum_{v \geq 0} (-1)^v \dim_k \left ( \Ext_{\gr A}^v(S_i,
S_j[u]) \right ) x^u.$$
\end{thm}

For $C_A(x)$ to be invertible as a matrix over $\Z[x]$, all entries
$c_{ij}$ above have to be polynomials (not proper power series), and
we get the following corollary.

\begin{cor}
$\det C_A(x)=1$ if and only if the entries $c_{ij}$ in $C_A^{-1}(x)$
are polynomials for all $1 \leq i,j \leq r$.
\end{cor}

If $A$ has finite global dimension, then $\det C_A(x)=1$. The
converse is not true, there are algebras $A$ of infinite global
dimension with $\det C_A(x)=1$.

\begin{example}\sloppy
Let $A$ be the path algebra $A=kQ/I$, where $Q$ is the quiver
$$\xymatrix{1\ar@/^/[r]^\alpha &
2  \ar@/^/[l]^\delta \ar@/^/[r]^\beta & 3\ar@/^/[l]^\gamma }$$ and
$I= \langle \rho \rangle$ is the ideal generated by the set of
relations $\rho =\{\beta \alpha, \gamma \beta, \beta \gamma, \delta
\gamma,\alpha \delta \alpha, \delta \alpha \delta\}$. Then $\det
C_A(x)=1$, but $A$ has infinite global dimension. Since $A$ is
monomial, Han's conjecture is known to hold for this algebra and
therefore $\hhdim A= \infty$.
\end{example}

For some classes of graded algebras, the combination $\det C_A=1$
and infinite global dimension is not possible. We can use Corollary
\ref{determinant} to prove that Han's conjecture holds for these
classes.

\subsection{Koszul algebras}

Let $A$ be a graded algebra with $A_0 \simeq k^{\times r}$. Such an
algebra is called \emph{Koszul} (\cite{Priddy}, \cite{BGS}) if
$\Ext_{\gr A}^v(S_i, S_j[u]) \neq 0$ implies $u=v$.

\begin{thm}\label{Koszul}
Let $A$ be a finite dimensional Koszul algebra over a field of
characteristic zero. If $A$ has infinite global dimension, then
$\hhdim A=\infty$.
\end{thm}

\begin{proof}
Suppose $\det C_A(x) =1$. Then $C_A(x)$ has an inverse $C_A^{-1}(x)$
whose entries are polynomials and \emph{not} proper power series.
Since our algebra is Koszul, we know that $\Ext_{\gr A}^v(S_i,
S_j[u])$ is nonzero only when $u=v$, and so we may simplify the
formula for the entries in the inverse matrix and get
$$c_{ij} = \sum_{u \geq 0} (-1)^u \dim_k \left ( \Ext_{\gr A}^u(S_i,
S_j[u]) \right ) x^u.$$ By assumption, all the entries in
$C_A^{-1}(x)$ are polynomials, hence all the higher graded extension
groups between simple $A$-modules vanish. Consequently, the global
dimension of $A$ is finite. This shows that the determinant of the
graded Cartan matrix of a Koszul algebra of infinite global
dimension is not one, and so by Corollary \ref{determinant} we are
done.
\end{proof}

\subsection{Local algebras}

If $A$ is graded finite dimensional with $A_0=k$, then the graded
Cartan determinant is by definition the same as the Hilbert
polynomial.

\begin{thm}
Suppose $k$ is of characteristic zero, and let $A$ be a graded
finite dimensional $k$-algebra with $A_0=k$. If $A$ has infinite
global dimension, then $\hhdim A=\infty$.
\end{thm}

\begin{proof}
If $\hhdim A < \infty$, then $\det C_A=1$ by Corollary
\ref{determinant}. If the Hilbert polynomial of $A$ is $1$, then
$A=k$ and $A$ has finite global dimension.
\end{proof}

This theorem can be seen as a special case of the following more
general result.

\begin{thm}\label{loop}
Suppose $k$ is of characteristic zero, let $Q=(Q_0,Q_1)$ be a finite
oriented quiver, and let $\mathtt J$ be the ideal in the path
algebra $kQ$ generated by the arrows. Furthermore, let $I$ be a
homogenous ideal in $kQ$ such that $\mathtt J^t \subseteq I
\subseteq \mathtt J^2$ for some $t$. If $Q$ contains a loop, then
$\hhdim kQ/I = \infty$.
\end{thm}

\begin{proof}
This follows from Corollary \ref{determinant} and the construction
of the graded Cartan matrix. The entries in $C_{kQ/I}(x)$ which are
not on the diagonal have constant term zero, so the only
contribution to the degree one coefficient in the determinant comes
from the product of the diagonal entries. The $i$th diagonal entry
is of the form $1+ n_i x + \cdots$, where $n_i$ is the number of
loops at vertex $i$ of $Q$. Therefore $$\det C_{kQ/I}(x)=1 +
\sum_{i=1}^{|Q_0|} n_i x + \text{higher terms}.$$ Consequently, if
$Q$ contains a loop, then $\det C_{kQ/I}(x) \neq 1$.
\end{proof}

\subsection{Cellular algebras}

Cellular algebras were introduced in \cite{Cell} (see also
\cite{KX1}) and are finite dimensional algebras which admit a
special kind of basis. A $k$-algebra $A$ is called a \emph{cellular
algebra} with cell datum $(\Lambda,M,C,i)$ if all of the following
three conditions are satisfied.
\begin{itemize}
\item[(C1)] The set $\Lambda$ is finite and partially ordered. Associated
with each $\lambda \in \Lambda$ there is a finite set $M(\lambda)$.
The algebra $A$ has a $k$-basis $\{C^{\lambda}_{S,T}\}$, where
$(S,T)$ runs through all elements of $M(\lambda) \times M(\lambda)$
for all $\lambda \in \Lambda$.
\item[(C2)] The map $i$ is a $k$-linear anti-automorphism of $A$
with $i^2=\mathtt 1_A$ which sends $C^{\lambda}_{S,T}$ to
$C^{\lambda}_{T,S}$.
\item[(C3)] For each $\lambda \in \Lambda$ and $S,T \in M(\lambda)$ and
each $a \in A$, the product $a \cdot C^{\lambda}_{S,T}$ can be
written as $$a \cdot C^{\lambda}_{S,T}=\sum_{U \in
M(\lambda)}h_a(U,S) \cdot C^{\lambda}_{U,T} + h',$$ where $h'$ is a
linear combination of basis elements with upper index $\mu$ strictly
smaller than $\lambda$, and where the coefficients $h_a(U,S)$ do not
depend on $T$.
\end{itemize}

Let $S_1, \ldots, S_r$ be a complete set of non-isomorphic simple
$A$-modules, and let $P_1, \ldots, P_r$ denote the corresponding
indecomposable projective modules. The \emph{(ungraded) Cartan
matrix} of $A$, denoted $U_A$, is the $r \times r$ matrix over $\Z$
where for all $1 \leq i,j \leq r$, the entry $(U_A)_{ij}$ is equal
to the composition multiplicity of $S_j$ in $P_i$. In \cite{KX2} we
find the following characterization of cellular algebras of finite
global dimension. For the definition of quasi-hereditary algebras,
see \cite{CPS}.

\begin{thm}\cite[Theorem 1.1]{KX2}
Let $A$ be a cellular algebra over a field. The following are
equivalent.
\begin{itemize}
\item[(a)] $A$ is quasi-hereditary,
\item[(b)] $A$ has finite global dimension,
\item[(c)] $\det U_A=1$.
\end{itemize}
\end{thm}

If $A$ is graded with $A_0 \simeq k^{\times r}$, then $U_A$ can be
obtained from $C_A(x)$ by evaluating for $x=1$. Therefore $\det U_A
=\det C_A(1)$.

\begin{thm}
Suppose $k$ is of characteristic zero, and let $A$ be a graded
cellular $k$-algebra such that $A_0$ is a product of copies of $k$.
If $A$ has infinite global dimension, then $\hhdim A=\infty$.
\end{thm}

\begin{proof}
If $\det U_A \neq 1$, then $\det C_A(x) \neq 1$, and so by Corollary
\ref{determinant} we have $\hhdim A=\infty$.
\end{proof}

\end{document}